\documentclass[10pt]{article}
\usepackage{amsmath,amssymb,amsthm,amscd}
\numberwithin{equation}{section}

\def\CC{{\mathbb C}}

\newtheorem{prop}{Proposition}[section]
\newtheorem{theo}[prop]{Theorem}

\newtheorem{cor}[prop]{Corollary}
\newtheorem{rem}[prop]{Remark}

\def\begeq{\begin{equation}}
\def\endeq{\end{equation}}

\def\and{\quad{\rm and}\quad}

\let\lra=\longrightarrow

\def\mapright\#1{\,\smash{\mathop{\lra}\limits^{\#1}}\,}

\begin {document}
\bibliographystyle{plain}
\title{Ricci flow on K\"ahler manifolds }
\author{X. X. Chen and G. Tian}
\date{July, 26, 2000}
 % Enter your date or \today between curly braces
\maketitle
\section{Introduction and main theorems}

In the last two decades, the Ricci flow, introduced by R.
Hamilton in \cite{Hamilton82}, has been a subject of intense
study. The Ricci flow provides an indispensable tool of deforming
Riemannian metrics towards canonical metrics, such as Einstein
ones. It is hoped that by deforming a metric to a canonical
metric, one can further understand geometric and topological
structures of underlying manifolds. For instance, it was proved
\cite{Hamilton82} that any closed 3-manifold of positive Ricci
curvature is diffeomorphic to a spherical space form. We refer
the readers to \cite{Hamilton93} for more information.

If the underlying manifold is a K\"ahler manifold, the normalized
Ricci flow in a canonical K\"ahler class \footnote{A K\"ahler
class is  canonical if the first Chern class is proportional to
this K\"ahler class. } preserves the K\"ahler class. It follows
that the Ricci flow can be reduced to a fully nonlinear parabolic
equation on almost pluri-subharmonic functions:
\[
\frac{{\partial \varphi }}{{\partial  t}} = \log \;\det \left(
{\frac{{\left( {\omega + \partial \bar \partial \varphi }
\right)^n }}{{\omega^n }}} \right) + \varphi  - h_{\omega},
\]
where $\varphi$ is the evolved K\"ahler potential; and  $\omega$
is the fixed K\"ahler metric in the canonical K\"ahler class,
while $Ric(\omega)$ is the corresponding Ricci form and
\[ Ric(\omega) - \omega =
\partial \bar
\partial h_{\omega}, \qquad {\rm and }\qquad \int_M (e^{h_{\omega} } - 1) \; {\omega}^n = 0.
\]

 Usually, this reduced flow is called the
K\"ahler Ricci flow. H.D. Cao \cite{Cao85} proved that the
K\"ahler Ricci flow always has a global solution. He also proved
that the solution converges to a K\"ahler-Einstein metric if the
first Chern class of the underlying K\"ahler manifold is zero or
negative. Consequently, he reproved the famous Calabi-Yau theorem
\cite{Yau78}. On the other hand, if the first Chern class of the
underlying K\"ahler manifold is positive, the solution of a
K\"ahler Ricci flow may not converge to any K\"ahler-Einstein
metric. This is because there are compact K\"ahler manifolds with
positive first Chern class which do not admit any
K\"ahler-Einstein metrics (cf.  \cite{futaki83}, \cite{Tian97}).
A natural and challenging problem is whether or not the K\"ahler
Ricci flow on a compact K\"ahler-Einstein manifold converges to a
K\"ahler-Einstein metric. It was proved by S. Bando
\cite{Bando84} for 3-dimensional K\"ahler manifolds and by N. Mok
\cite{Mok88} for higher dimensional K\"ahler manifolds that the
positivity of bisectional curvature is preserved under the
K\"ahler Ricci flow. A long standing problem in the study of the
Ricci flow is whether or not the K\"ahler Ricci flow converges to
a K\"ahler-Einstein metric if the initial metric has positive
bisectional curvature? In view of the solution of the Frankel
conjecture by S. Mori \cite{Mori79} and Siu-Yau \cite{Siuy80}, we
suffice to study this problem on a K\"ahler manifold which is
biholomorphic to $\CC P^n$. Since $\textsc{C} P^n$ admits a
K\"ahler-Einstein metric, the above problem can be restated as
follows: On a compact K\"ahler-Einstein manifold, does the
K\"ahler Ricci flow converge to a K\"ahler-Einstein metric? In
this note, we announce an affirmative solution to this problem.

\begin{theo}\cite{chentian001}\cite{chentian002}  Let $M$ be a K\"ahler-Einstein
manifold with positive scalar curvature. If the initial metric has
nonnegative bisectional curvature and positive at least at one
point, then the K\"ahler Ricci flow will converge exponentially
fast to a K\"ahler-Einstein metric with constant bisectional
curvature.
\end{theo}

\begin{rem} The above theorem in complex dimension $1$ was proved
first by Hamilton \cite{Hamilton88}.   B. Chow \cite{Chow91}
later showed that the assumption that the initial metric has
positive curvature in $S^2$ can be removed since the scalar
curvature will become positive after finite time anyway.
\end{rem}

\begin{cor}
The space of K\"ahler metrics with non-negative bisectional
curvature is path-connected. Moreover, the space of metrics with
non-negative curvature operator\footnote{In \cite{Hamilton86},
Hamilton proved that the positivity of the curvature operator is
preserved under Ricci flow in any compact manifold.} is also
path-connected.
\end{cor}

\begin{rem}
Using the same arguments, we can also prove the version of our
main theorem for K\"ahler orbifolds.
\end{rem}

\begin{rem}
What we really need is that the Ricci curvature is positive.
Since the condition on Ricci may not be preserved under the Ricci
flow, in order to have the positivity of the Ricci curvature, we
will use the fact that the positivity of the bisectional
curvature is preserved.
\end{rem}

\begin{rem} We need to assume the existence of K\"ahler-Einstein
metric because we will use a nonlinear inequality from
\cite{tian98}. Such an inequality is nothing but the
Moser-Trudinger-Onofri type inequality if the K\"ahler-Einstein
manifold is the Riemann sphere.
\end{rem}

\section{Outline of Proof}

 The standard methods in the Ricci flow involve the
pointwise estimate of curvature by using its evolution equation,
the blow-up analysis and the classification of possible
singularity type. Unfortunately, there are only a few examples of
proving convergence of the Ricci flow to a canonical one (cf.
\cite{Hamilton82} \cite{Hamilton86} and \cite{Hu85} et al).  One
of the main reasons is that it is very hard to detect geometric
information from a singular model arisen from a blow-up
analysis.\\

 In order to overcome this difficulty, we define a set of new
functionals $E_k(\varphi)\;(k=0,1,2,\cdots n).\;$   The leading
term of $E_k\;$ ($\; k = 0,1\cdots n$) is
   \[
   \displaystyle \int_M \; \left(\ln \det \left( {{\omega_{\varphi}}^n \over {\omega}^n} \right) -
   h_{\omega}\right) \displaystyle \sum_{i=0}^k Ric(\omega_{\varphi})^i \wedge {\omega}^{k-i}
   \wedge {\omega_{\varphi}}^{n-k}.
   \]
Here $\omega_{\varphi} = \omega +  \partial \bar \partial
\varphi  $ is the K\"ahler metric determined by the K\"ahler
potential $\varphi.\;$ The Euler-Lagrange equation for the
functional $E_k$($\; k = 0,1\cdots n$) is
   \[
     \triangle_{\omega_{\varphi}} \; \left( {{Ric(\omega_{\varphi})^k  \wedge {\omega_{\varphi}}^{n-k}} \over {\omega_{\varphi}}^n}\right)
   - {{n-k}\over {k+1}} \left( {{Ric(\omega_{\varphi})^{k+1}  \wedge {\omega_{\varphi}}^{n-k-1}} \over {\omega_{\varphi}}^n}
   \right) = - {{n-k}\over {k+1}}  \displaystyle \int_M {\omega_{\varphi}}^n.
   \]
   If $k=0, \;$ this leads to the usual constant scalar curvature metric equation.  In general,
   critical metrics may not  have constant scalar curvature.\\

  The derivative of the functional $E_k$  over a one
   parameter family of automorphisms gives rise to a holomorphic
   invariant ${\cal F}_k\;( k=0,1,2,\cdots n).\; $
   One shall note that these functionals $E_k$ and holomorphic invariants  ${\cal F}_k $
   are defined in any K\"ahler class
   of any K\"ahler manifold. However, if the K\"ahler class is  canonical on a K\"ahler-Einstein manifold,
   then these invariants vanish simultaneously. It follows that on a K\"ahler-Einstein manifold,
   these functionals are invariant under action of  automorphisms.  This important
property enables us to modify the
   flow by automorphisms so that we can apply the fully nonlinear inequality of Tian \cite{tian98}. In fact, this inequality
   is  a generalized version of the Moser-Trudinger-Onfri  inequality\footnote{In $S^2$, we need
to require a function to be perpendicular to
   the first eigenspace of the standard metric in $S^2$ in order to get better estimates. We need to do exactly
   the same here: we need to adjust the K\"ahler Ricci flow by automorphisms so that the evolved K\"ahler potential is
   always perpendicular to the first eigenspace of a fixed K\"ahler-Einstein metric.}. Using the inequality,
   we can derive a uniform low bound of the evolved volume form along the modified K\"ahler Ricci flow;
   Consequently, we can derive a uniform
   lower bound of these functionals along the
   K\"ahler Ricci flow.  \\

   Furthermore, these
   functionals essentially decrease along the K\"ahler Ricci flow. Since the functional $E_k$ has
    a uniform lower bound over the entire K\"ahler Ricci flow, by computing its derivative, we
   have
   \[
   \int_{t=0}^{\infty} \; \int_M \; R(\omega_{\varphi}) ({Ric(\omega_{\varphi})}^k - \omega^k) \wedge
   \omega^{n-k}\; d\,t< C
   \]
for some uniform constant $C.\;$
   In particular when $k=1,\; $ we obtain
   \[
 \int_{t=0}^{\infty} \; \int_M \; (R(\omega_{\varphi})-r)^2 {\omega_{\varphi}}^n \; d\,t <
 C,
\]
where  $r$ is the average of the scalar curvature --- a constant
depending only on the K\"ahler class. \\

  This means that $\int_M \; (R(\omega_{\varphi})-r)^2 {\omega_{\varphi}}^n $ is  small for almost all the time. In complex
  dimension 2,   combining this integral estimate with Cao's Harnack inequality for the scalar curvature
    and a generalized version of Klingenberg  estimate on injective radius, we can prove that the scalar curvature
    is uniformly bounded from above over the entire K\"ahler Ricci flow.
    Once the scalar curvature is uniformly bounded,
    we then follow the standard arguments in the theory of parabolic equations to
    deduce the exponential convergence to a K\"ahler-Einstein metric for this flow. Therefore, we prove
Theorem 1 in complex dimension 2.  For high
    dimensional K\"ahler manifolds, it is more complicated to  control the scalar curvature and
     we need to use some new ingredients from the theory of pluri-subharmonic functions. \\

    The detailed proof of these results will appear elsewhere.
%\bibliography{test}

\noindent Department of mathematics, Princeton University,
Princeton, NJ 08544, USA; \\
\noindent xiu@math.princeton.edu\\

\noindent Department of mathematics, MIT, Cambridge, MA
02139-4307, USA;\\
\noindent tian@math.mit.edu

 \end{document}